\documentclass[11pt]{article}
\usepackage{graphicx}
\usepackage{wrapfig}
\usepackage{manyfoot}

\usepackage{amsmath}
\usepackage{url}

\newcommand\sLP{\\[\smallskipamount]}
\newcommand\mLP{\\[\medskipamount]}
\newcommand\mPP{\\[\medskipamount]\indent}
\newcommand\bLP{\\[\bigskipamount]}
\newcommand\al\alpha
\newcommand\be\beta
\newcommand\Ga{\Gamma}
\newcommand\De{\Delta}
\newcommand{\hyp}[5]{\,\mbox{}_{#1}F_{#2}\!\left(
  \genfrac{}{}{0pt}{}{#3}{#4};#5\right)}

\numberwithin{equation}{section}

\DeclareNewFootnote{A}[arabic]

\begin{document}
\title{Dick and Liz Askey's visit to U.S.S.R.\ in 1987, and
how the discrete Askey scheme also originated in Russia}
\author{Tom H. Koornwinder}
\date{}
\maketitle
\begin{abstract}
This paper describes how the discrete Askey scheme independently arose
in Russia and how Askey learned about this. In particular, Askey met main
characters in this story, namely Gel'fand and Suslov as well as Nikiforov
and Uvarov, during his trip to U.S.S.R.\ in September 1987. The paper
describes this trip in some detail, in particular based on the diary of
Askey's wife Liz, who accompanied him.
Dick and Liz Askey continued their trip by visits to Japan,
Australia and India. Schedules of these three visits are also given.
\end{abstract}

\section{Introduction}
During September 1987 -- January 1988 Dick Askey and his wife Liz
made a trip
successively to the U.S.S.R., Japan, Australia and India.
After the visit to the
U.S.S.R. they briefly returned to their home in Madison, but the other three
countries were visited without interruption.
Liz Askey wrote four travel diaries,
one for each country. These were originally handwritten. Afterwards a few
typewritten copies were circulated, also to me. Recently
I have made available edited and annotated versions of these diaries on the
web\footnoteA{\url{https://staff.fnwi.uva.nl/t.h.koornwinder/specfun/#DickAskey}}.

The diaries do not contain any mathematics, but Liz describes encounters with
(sometimes famous) mathematicians at which she was present.
The main focus in the diaries is on the cultural activities
during the trip, sometimes by Dick and Liz
together and often by Liz alone, while Dick discussed math.

Their travel schedule was as follows.
\begin{itemize}
\item
September 1--21: U.S.S.R.\ (Moscow and Leningrad)
\item
back in Madison for two weeks
\item
October 6 -- November 14: Japan
\item
November 15 -- December 10: Australia (Sydney and Canberra)
\item
December 11 -- January 8: India (Ramanujan birthday centennial)
\item
January 9--15: London, UK
\end{itemize}
In this paper I first give some context, including a short biography of
Liz Askey, some words about Dick Askey's mathematical status when he embarked
on this trip, and an account of
important mathematical developments in those years which touched his own area
of interest. Then follows in Section 3 a description of work done in Russia
on elements of the discrete Askey scheme, either original or independent
from work done elsewhere. This ranges from Chebyshev in the 19th century until
Sergei Suslov and his collaborators in the 1980s. The aim of this section
is to put the next section on Askey's trip to U.S.S.R./ in a better
perspective. For both sections Sergei Suslov's eyewitness communications
to me were of inestimable value.
See also his recollections \cite{Suslov25} of Askey

The paper concludes with three short sections on the Askey trips to
Japan, Australia and India, only listing the places visited
and the conferences attended.
For more details on these trips see the
online edited diaries by Liz Askey, the slides\footnoteA%
{\url{https://staff.fnwi.uva.nl/t.h.koornwinder/art/sheets/2022_OPSFA.pdf}}
of my lecture on June 15, 2022 at OPSFA 16,
and about India the eyewitness account
by George Andrews \cite{Andrews25}.
\section{Some context}
\subsection{Liz Askey}
Elizabeth Hill ``Liz'' Askey\footnoteA{\url{https://www.legacy.com/us/obituaries/madison/name/elizabeth-askey-obituary?id=33536679}}
was born November 18, 1936 in Aberdeen, Washington.
She got married to Dick Askey on June 14, 1958. She died January 29, 2022,
surviving Dick for 2$\frac14$ years. They had two children: a son and a
daughter. Liz was an undergraduate at Bryn Mawr College majoring in history.
She met Dick at a mixer at her College with Princeton University.

After her children were in school, Liz studied part time at the University of
Wisconsin-Madison and received a master's degree in urban and regional planning.
She went on to receive another master's degree from the School of Library and
Information Studies. Children's literature always had her special interest. She was
a good writer and storyteller.

Liz often accompanied Dick at his travels. Many readers
will have met Liz at a conference where Dick was speaking.
\subsection{Dick Askey's mathematical achievements and interests by 1987}
By the time of his 1987 trip Dick Askey was at the height of his
mathematical power. In 1984 his AMS Memoir \cite{Askey-Ismail84} with Mourad
Ismail had appeared,
and in 1985 his AMS Memoir~\cite{Askey-Wilson85} with Jim Wilson,
introducing the Askey--Wilson polynomials.
This 1985 Memoir also introduced the celebrated \emph{Askey scheme} (yet
with only a 2-parameter class of continuous Hahn polynomials).
As I mentioned in \cite[Foreword, p.~vi]{Koekoek2010}, 
Askey got the idea of such a scheme from an Oberwolfach lecture in 1977
by Michael Hoare, who presented a chart
from his 1977 paper \cite[Fig.~2]{Cooper-Hoare-Rahman77}
with Cooper and Rahman which goes down from
the Hahn polynomials to the Hermite polynomials.

A $q$-analogue of the Askey scheme was essentially described in the
Askey--Wilson Memoir, but it was not yet displayed as a graph.
In fact, Askey then knew already for a long time, from the Thesis
\cite[Chapter~III]{Wilson78} in 1978 by his student Jim Wilson, that the
($q$-)Askey scheme did not yet reach the top: Wilson had shown that there are
lying biorthogonal rational functions above the Racah polynomials.
It lasted until 1991 before
Wilson published this in a journal \cite{Wilson91}. Until then
Askey did not speak widely about the result.

The ($q$-)Askey scheme could also be extended by considering the
corresponding \emph{associated} orthogonal families. This would culminate
in 1991 when Ismail and Rahman \cite{Ismail-Rahman91}
presented the associated Askey--Wilson polynomials.

In the same year 1985 as the Askey--Wilson Memoir, Louis
de Branges' proof~\cite{deBranges85} of the Bieberbach conjecture was
published. This used
on p.150 the Askey--Gasper inequality \cite[(1.17), (3.1)]{Askey-Gasper76} as a
crucial step.

Titles of Askey's lectures in the first three visited countries
have not been preserved, so we can only speculate about the topics.
I guess he gave quite a few lectures on $q$. Also Ramanujan's work and legacy
may have been among the topics. This would have been a good preparation
for the Ramanujan centennial
in India, December 1987, where he had to speak at least five times at
memorial conferences.
\subsection{Some new mathematical developments by 1987}
Just before or right in 1987 some new theories came up which were also relevant for
orthogonal polynomials and special functions:
\begin{itemize}
\item
Quantum groups (Drinfel'd \cite{Drinfeld87}, ICM Berkeley, 1986)
\item
Heckman--Opdam polynomials \cite{Heckman-Opdam87} (1987, 1988)
\item
Macdonald polynomials (manuscripts circulated in 1987)
\item
Dunkl operators \cite{Dunkl89} (1989, submitted in 1987)
\item
Gelfand hypergeometric functions \cite{Gelfand86} (1986)
\item
Wavelets (ondelettes \cite{Lemarie-Meyer86}, 1986)
\end{itemize}
\noindent
While embarking on his trip, Dick may not yet have been fully aware of all
these, but he will certainly have learned more about these topics from people
he met during the trip.
\section{The discrete Askey scheme in Russia}
\label{S3}
\subsection{Chebyshev's 1875 paper}
For about a century Chebyshev's paper  \cite{Chebyshev1875} (1875) was
mainly remembered because of the \emph{discrete Chebyshev polynomials}
\cite[\S2.8]{Szego67}, \cite[\S10.23]{HTF2}, i.e.,
the orthogonal polynomials for constant weights on the finite set
$\{0,1,\ldots,N\}$. In fact, as described by Roy \cite{Roy93},
Chebyshev wrote several earlier papers \cite{Chebyshev1858} (1858),
\cite{Chebyshev1859} (1859), \cite{Chebyshev1864} (1864)
on interpolation\footnote{See also on R. J. Pulskamp's webpage\\
\url{https://probabilityandfinance.com/pulskamp/Chebyshev/Chebyshev.html}
a listing of such papers
under the heading ``Interpolation by Method of Least Squares'',
with links to English translations.},
where he discusses these polynomials.
As described in \cite[\S3]{Roy93}, Chebyshev obtains his polynomials
as the denominator polynomials in the continued fraction expansion
\cite[Chapter III, \S4]{Chihara78}
of
\[
\frac1{x-1}+\frac1{x-2}+\cdots+\frac1{x-m+1}\,.
\]
He shows the orthogonality and he derives a finite difference variant of
the Rodrigues formula for Legendre polynomials.

However, the two-parameter generalization of these orthogonal polynomials,
written in current notation  \cite[\S9.5]{Koekoek2010} as
\begin{equation}
Q_n(x;\al,\be,N):=\hyp32{-n,n+\al+\be+1,-x}{\al+1,-N}1,\quad n=0,1,\ldots,N,
\end{equation}
is introduced in the later sections 9 and 10 of \cite{Chebyshev1875}.
There Chebyshev obtains polynomials
$\phi_n(x)$ as orthogonal polynomials for weights
\[
\theta^2(x):=\frac{\Ga(x+\al)}{\Ga(x)}\,\frac{\Ga(m-x+\be)}{\Ga(m-x)}=
\Ga(\al+1)\Ga(\be+1)\,\frac{(\al+1)_{x-1}}{(x-1)!}\,
\frac{(\be+1)_{m-x-1}}{(m-x-1)!}
\]
on the finite set $\{1,2,\ldots,m-1\}$. He shows that these OPs 
can be given by the formula
\[
\phi_n(x)=\De^n\left(\frac{\Ga(x+\al)}{\Ga(x-n)}\,
\frac{\Ga(m-x+\be+n)}{\Ga(m-x)}\right)\Big/
\left(\frac{\Ga(x+\al)}{\Ga(x)}\,\frac{\Ga(m-x+\be)}{\Ga(m-x)}\right).
\]
(what we now call a \emph{Rodrigues formula})
because the right-hand side defines polynomials of degree $n$ satisfying
the same orthogonality relations as $\phi_n(x)$.
He continues by giving an evaluation of the quadratic norm of $\phi_n(x)$ and
he concludes by deriving an expression for $\phi_n(x)$ as a finite sum
which we can write as
\begin{align*}
\phi_n(x)&=(\al+1)_n (m-n-x)_n\,\hyp32{-n,-n-\be,-x+1}{\al+1,m-n-x}1\\
&=(m-n-1)_n (\al+1)_n\,\hyp32{-n,n+\al+\be+1,-x+1}{\al+1,-m+2}1\\
&=(m-n-1)_n (\al+1)_n\,Q_n(x-1;\al,\be,m-2).
\end{align*}
Here the second equality comes from a transformation formula for
terminating ${}_3F_2(1)$ which is a special case of \cite[3.8(1)]{Bailey35}:
\begin{equation}
\hyp32{-n,b,c}{d,e}1=\frac{(e-b)_n}{(e)_n}\,\hyp32{-n,b,d-c}{d,b-e-n+1}1.
\end{equation}
Read more about Chebyshev's work in \cite{Butzer99}.

The reason that the case of general $\al,\be$ in \cite{Chebyshev1875} was
so long ignored is probably that the interest in these polynomials was mainly
for applications in interpolation and numerical analysis, where the case
$\al=\be=0$ is in particular useful.

The case of general $\al,\be$ was rediscovered by Hahn \cite{Hahn49}, but
mostly in the $q$-setting. His type~I gives rise to what we now call
$q$-Hahn polynomials, and their limit case for $q\to1$ yields the Hahn
polynomials, only briefly treated at the end of Hahn's paper.
They are also briefly given in \cite[\S10.23]{HTF2}, and they are treated
in more detail in \cite{Weber52} and \cite{Karlin61}.
All these papers also mention Bateman's polynomials
\cite[(18), (19)]{Bateman34} (also considered by Hardy) as a special case of
the Hahn polynomials,
but these are in fact special continuous Hahn polynomials, which only
formally can be seen as Hahn polynomials \cite[\S9.4]{Koekoek2010}.
The Bateman polynomials were extended by Pasternack.
See \cite{Koelink93} and \cite[pp.~8, 9]{Suslov-Trey08}
for history and further references.

It was probably first in 1978 that Askey \cite[p.8]{Mackie78} gave the
credit for the general Hahn polynomials to Chebyshev. He gives a longer
discussion of Chebyshev's paper in \cite[p.866]{Askey82}.
\subsection{The discrete Askey scheme and its relation with group theoretical
physics}
Wilson obtained in his Thesis \cite{Wilson78} (1978) new discrete orthogonal
polynomials with respect to weights on a finite set which realized the
known orthogonality of the $6j$ symbols (or Racah coefficients) for the group
SU(2).
Thus the \emph{Racah polynomials} were born, on top of the discrete
Askey scheme, and with Hahn polynomials as an immediate limit case. He also
gave the version having continuous orthogonality, which was later called after
him: the \emph{Wilson polynomials}. He announced these results in
\cite{Wilson77} (1977) and he published them in a journal \cite{Wilson80}
in 1980. Already in 1979 he published together with Askey \cite{Askey-Wilson79}
a $q$-version of the Racah polynomials.

Around that time, in 1978, Askey \cite{Askey-Wilson79} was aware that the
$3j$ or Clebsch--Gordan coefficients for SU(2) could be expressed in
terms of Hahn polynomials, but the details were not published anywhere.
Then the author \cite{Koornwinder81} (1981) published this.
Also, in a publication \cite{Koornwinder82} one year later he observed
that the matrix elements of irreducible unitary representations of SU(2)
could not only be expressed in terms of Jacobi polynomials (a result due to
Wigner), but also in terms of Krawtchouk polynomials, realizing the
column orthogonality.
\subsection{Meanwhile in Moscow}
Gel'fand, Minlos \& \v{S}apiro \cite[Supplement III]{Gelfand58} observed in
1958 that Clebsch--Gordan coefficients could be considered as discrete
analogues of Jacobi polynomials because they can be expressed by a Rodrigues
type formula like the one for the Jacobi polynomials, but with the derivative
replaced by the difference operator. However, they did not express the
CG coefficients in terms of Hahn polynomials; they were unaware of these
polynomials. According to Graev (in a conversation with Suslov),
it was really Gel'fand who obtained this formula. Gel'fand was very excited
about this result.
He encouraged people around him to do further work on it, but nobody
was interested.
This disinterest was also apparent in a meeting of the Gelfand seminar
much later, where the lecture was on orthogonal polynomials.
As mentioned by Retakh \cite[p.164]{Retakh13},
only the senior people in the audience paid attention to the speaker.

What was left open in \cite{Gelfand58} was taken up in 1982
at the Kurchatov Institute by Smorodinskii and
his student Sergei Suslov \cite{Smorodinskii82a}: they expressed the
CG coefficients in terms of Hahn polynomials, independently from
\cite{Koornwinder81}. In the same year they rediscovered \cite{Smorodinskii82b}
the Racah polynomials and their connection with the $6j$ symbols,
independently from Wilson's work. Somehow Klimyk had pointed the authors to the
preprint version of \cite{Koornwinder81}, a report of the Mathematical
Centre in Amsterdam which was also available in Moscow.
So they learned about the results in \cite{Koornwinder81} and, via a reference
there, also about Wilson's paper \cite{Wilson80}. In a supplement to
\cite{Smorodinskii82b} they acknowledged these results.
\subsection{The evolution of the book by Nikoforov and Uvarov}
In the 1970s and 1980s Nikiforov and Uvarov, working at the Keldysh
Institute\footnoteA{\url{https://en.wikipedia.org/wiki/Keldysh_Institute_of_Applied_Mathematics}}
in Moscow, wrote two successive books \cite{N-S74}, \cite{N-S78} on
special functions, which went through various editions:
\mLP
1974: \emph{Foundations of the theory of special functions} (in Russian),\\
1977: French translation,\\
1978: \emph{Special functions of mathematical physics} (in Russian),\\
1983: French translation,\\
1984: second Russian edition,\\
1988: English translation.
\mLP
The first book \cite{N-S74} already had material (Section 17) on classical
discrete orthogonal polynomials, notably Hahn polynomials and their
Rodrigues formula, but nothing about their interpretation as
Clebsch--Gordan coefficients. The same is true for the first edition of
their second book \cite{N-S78}. 

Around 1982 Suslov got in touch with the 
authors (they had to meet outside the Keldysh Institute, because security
protocol did not allow non-employees to enter the building). Together
they published several reports in the Keldysh Institute preprint series.
Notably \cite{N-S82} gives again the results of \cite{Smorodinskii82a}, with
due references to \cite{Smorodinskii82a} and \cite{Koornwinder81},
and \cite{N-S-U82} deals with Racah polynomials (but no $6j$ symbols)
with due references to \cite{Smorodinskii82b}, \cite{Askey-Wilson79} and
\cite{Wilson80}.
Two more related papers by Suslov alone appeared in 1983 in
Soviet J.~Nuclear Phys.
No wonder that these results were also included in the next editions of
\cite{N-S78}. The 1983 French translation had the Hahn polynomial
interpretation as CG coefficients, and it treated the Racah polynomials
The 1984
second Russian edition also had the Racah polynomial interpretation as
$6j$ symbols. However, the just mentioned references in \cite{N-S82} and
\cite{N-S-U82} were absent in these two editions of \cite{N-S78}, but
\cite{Askey-Wilson79}, \cite{Wilson80} and \cite{Koornwinder81}
occurred as references in the 1988 English translation. This edition has also
$q$-analogues of the discrete classical OPs.
Suslov, together with a few
others, was thanked in the Introduction ``for helpful comments on the contents
of the book''.
\subsection{The ICM 1983 in Warsaw}
Probably Askey's first contacts with Soviet mathematicians occurred at the
ICM in Warsaw \cite{ProcICM84},
August 16--24, 1983, where Dick was an invited speaker.
He met there with Sergey Nikolsky (chair of the Soviet delegation to the ICM)
and Andrey Gonchar. Both were members of the Steklov
Institute\footnoteA{\url{https://en.wikipedia.org/wiki/Steklov_Institute_of_Mathematics}}
of the Soviet Academy of Sciences in Moscow.
Two other Moscow mathematicians which Askey had liked to meet,
Arnold Nikiforov and Israel Gel'fand,  were
not allowed to go there, because they were involved with work done at the
Keldysh Institute
in Moscow for the Soviet atomic bomb
project. However, Nikiforov had already submitted the title of his short
communication, and his lecture was scheduled. Then Dick seems to have
done his presentation.

It is not clear to which extent Askey was then, in 1983, familiar with the
work of Nikiforov, since the latter had until then only published in Russian.
But certainly, if they would have met then in Warsaw, there would have been
a lot to discuss.
\subsection{Suslov's collaboration with Atakishiyev: continuous
Hahn polynomials}
While Nikiforov and Uvarov were focused on discrete orthogonal polynomials,
ignoring the continuous families,
Suslov had also started a collaboration with Natig Atalishiyev
from Baku, then within the U.S.S.R. They had an open eye for developments
in the continuous case. An impactful result of their joint
work had to do with a paper \cite{Askey-Wilson82} (1982)
by Askey and Wilson, who gave a 2-parameter
family of continuous Hahn polynomials as a limit case of the
Wilson polynomials, which Wilson \cite{Wilson80} had introduced in 1980.
Their 2-parameter family was incorporated as part of the Askey scheme in
\cite{Askey-Wilson85} (1985). Atakishiyev and Suslov \cite{A-S85} pointed
out in the same year that the continuous Hahn polynomials could in fact
be extended
to a 3-parameter family, which was confirmed by Askey \cite{Askey85}
in a letter to the same journal, and which from then on was part of the
Askey scheme.
\subsection{A successor to \cite{N-S78}:
the book by Nikiforov, Suslov and Uvarov}
In 1983 Nikiforov and Uvarov also started working with Suslov on a book
focusing on discrete classical orthogonal polynomials including their
$q$-analogues. This book \cite{N-S-U85b} would appear in 1985.
Elements of the draft were already included in the later editions of
\cite{N-S78} (see above). Around 1984 Nikiforov was invited by Gel'fand
to tell about this work in his seminar, in a 20 minutes lecture after the
main speaker. Gel'fand showed much interest and also discussed things there
with Suslov. On the invitation of Gel'fand the paper \cite{N-S-U85a}
by Nikiforov, Suslov and Uvarov appeared in his journal.
At the end of the paper the authors thank Gel'fand and the participants in
the seminar he directs for helpful discussion.

I only have access to the English translation of the book \cite{N-S-U85b},
which appeared
in 1991. A short paper by Nikoforov and Suslov with some of the main results
from the book appeared in English \cite{N-S86} in 1986.
Askey, in his review MR0824673 of this paper in \emph{Mathematical Reviews},
compared their approach with the approach by Askey and Wilson.
The Russian approach was motivated by physics and
emphasized the ($q$-)difference operator going with such
polynomials, as an alternative to the American approach using
the three-term recurrence relation coming from the \mbox{($q$-)}hypergeometric
function. Askey continues:
``Most of their work was done independently, for while they knew of the first
Askey--Wilson paper \cite{Askey-Wilson79},
we did not include the divided difference form of the
Rodrigues formula there.'' But according to Suslov (personal
communication) they had not studied \cite{Askey-Wilson79} closely,
although they had included it in \cite{N-S86} as a reference.

The English translation of \cite{N-S-U85b} was certainly up-to-date with
references
to work done elsewhere. It had a very good reception and it is still considered
as a standard reference. The Preface gives credits to Suslov for having written
specific parts of the book, but in fact he contributed much more than
is mentioned there.
The Foreword to the Russian edition, also included in the translation,
was written by M.~I.~Graev in 1984. He explicitly gives credit to Gel'fand
for having noted in \cite{Gelfand58} the deep analogy between classical
orthogonal polynomials of continuous and discrete arguments, also in
connection with the representations of the rotation group. The
Acknowledgments express gratitude to Academicians S.~M.~Nikolsky and
A.~A.~Gonchar, ``who took part in the discussions of the authors' report as
well as the report of R. Askey in the Steklov Mathematical Institute''.
As will be described in the next section, these discussions took place
in September 1987, when the Russian edition had already appeared.
\section{U.S.S.R., September 1--21, 1987}
\begin{minipage}{11cm}
\paragraph{Overview}
The Askeys left August 31 from Madison (see photo) for a three weeks stay
in the U.S.S.R. After arrival in Moscow they stayed there only one night. and
then took the night train to Leningrad, where they stayed during
September \mbox{3--9}. After return to Moscow.
they stayed at Hotel Akademicheskaya, likely arranged by
the Moscow Steklov Institute. In Leningrad they stayed at Hotel Moskva,
arranged by Sergey Khrushchev of the Leningrad Steklov Institute.
\end{minipage}
\qquad\qquad
\begin{minipage}{2.5cm}
\includegraphics[width=2.5cm]{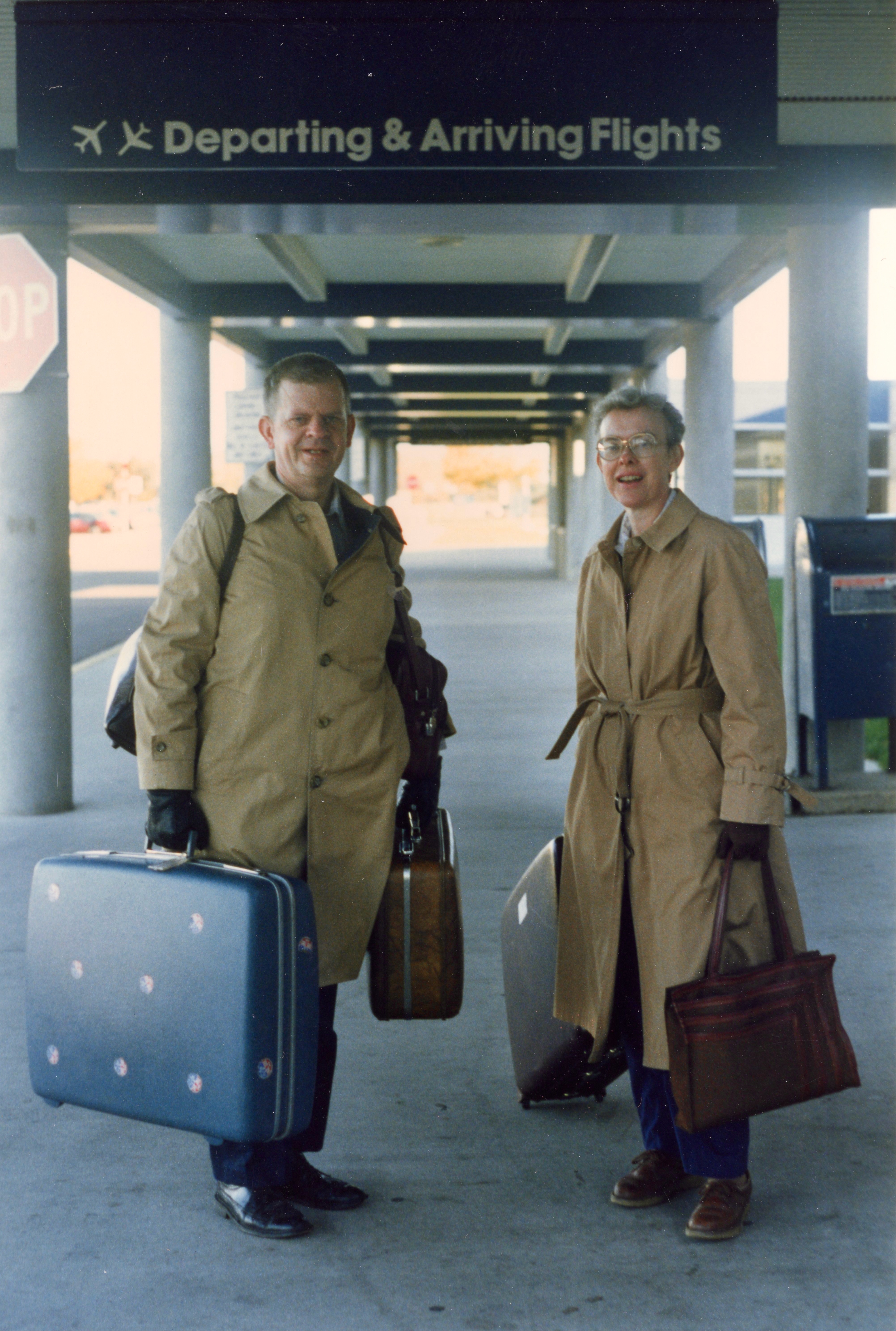}
\end{minipage}
\paragraph{September 1: Arrival at the airport in Moscow}\quad\mLP
\begin{minipage}{8.0cm}
\textsl{We arrived at the Moscow airport on September 1, after the standard
three hours of sleep on the plane, and inched our way through customs.
The man ahead of us, who was evidently visiting Russian relatives and who spoke
Russian himself, seemed to have brought in (as Dick phrased it) the contents of
a K-Mart store.
The customs man removed a tall pile of items including videotapes and a boom
box and the visitor tried his best to argue him out of all those questionable
items.
Finally the discussion, the man, and the stack of ``contraband'' were taken to
some quieter corner.
We saw the man rejoin his family later but we never did find out whether he
argued Russian customs out of anything. Having whetted his appetite on the
man ahead of us, the customs
official was ready for us; he went through everything with the proverbial
``fine toothed comb''.}
\end{minipage}
\quad
\begin{minipage}{7.6cm}
\includegraphics[width=7.6cm]{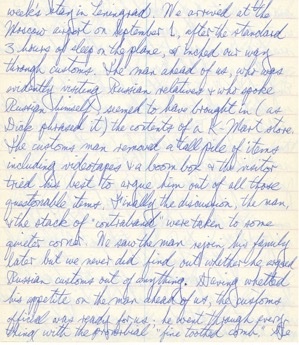}
\end{minipage}
\\
\begin{minipage}{8.5cm}
\textsl{He wondered about the hearing aid which Dick had brought for Gel'fand,
but the sticky point proved to be the anthology of twentieth century
Russian literature --- specifically a story by Nabokov.
A supervisor stood reading for some time before letting us take the book
into the country.}
\end{minipage}
\quad
\begin{minipage}{7.6cm}
\includegraphics[width=7.6cm]{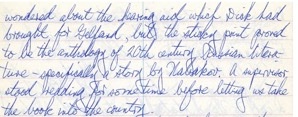}
\end{minipage}
\mLP
\begin{minipage}{7.8cm}
\paragraph{September 2: Visit to Gel'fand's home}\quad\sLP
\textsl{Our guide from the night before found us a taxi and took us to the
apartment building in which Gel'fand and his family lived.
We were only going to stay fifteen minutes since Gel'fand's wife was ill.
He was delighted to see Dick and the conversation went on much longer than 15
minutes.
He was interrupted constantly to answer the telephone and we eventually
learned that it was his 74th birthday.}
\end{minipage}
\qquad
\begin{minipage}{7cm}
\includegraphics[width=3.5cm]{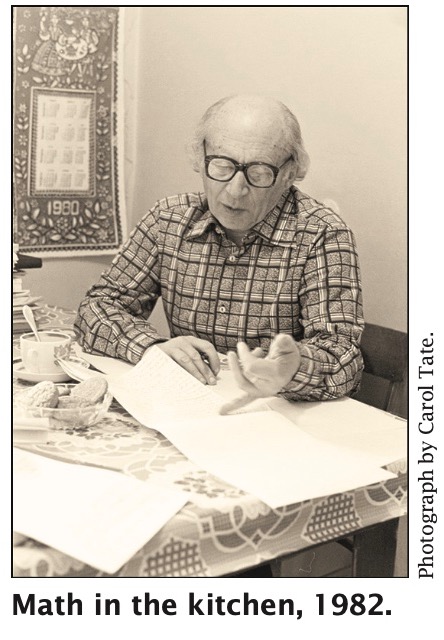}
\end{minipage}
\mPP
Israel Gel'fand%
\footnoteA{\url{https://en.wikipedia.org/wiki/Israel_Gelfand},\\
\url{https://mathscinet.ams.org/mathscinet/MRAuthorID/189130}} \cite{Retakh13}
(1913--2009) was one of the most prominent post-war
mathematicians in Moscow, active in many fields, notably group theory,
representation theory and functional analysis.
His work \cite{GelfandLevitan55} with Levitan in 1955, combined with a
publication by Marchenko in the same year, gave rise to the
Gel'fand--Levitan--Marchenko
integral equation, which lies on the basis of inverse scattering theory and is
important for soliton theory.
In 1987 Gel'fand had just started
work, together with others in his group, on the multivariable hypergeometric
functions named after him. The Gelfand seminar on Mondays was famous, but
frightening for some speakers and participants.
In 1989 he emigrated to the United States.
\begin{wrapfigure}{R}{0.5\textwidth}
\includegraphics[width=0.5\textwidth]{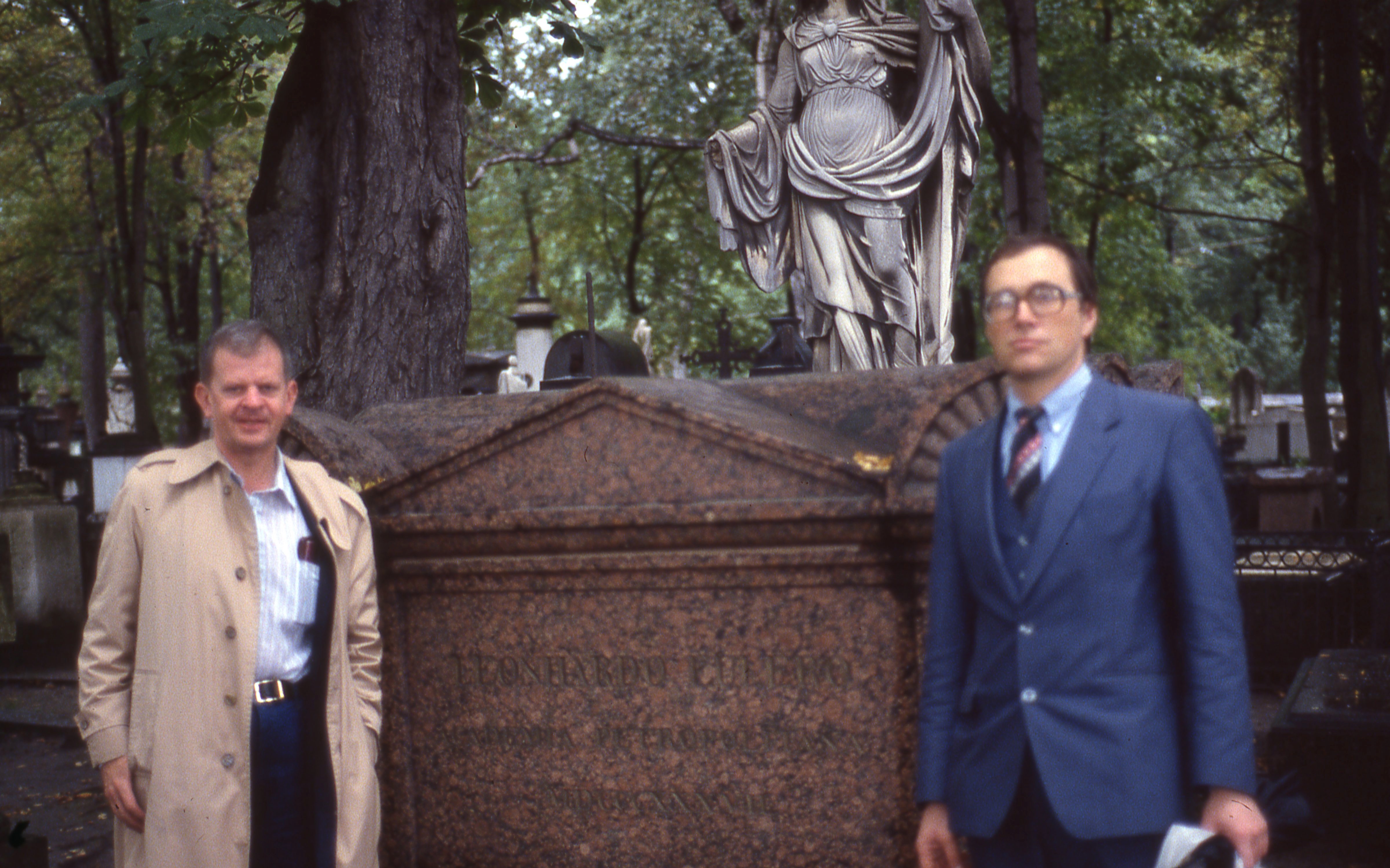}
\end{wrapfigure}
\paragraph{September 3--9: Visit to Leningrad}\quad\sLP
The Moskva Hotel in Leningrad is on the eastern end of Nevsky Prospekt.
In the immediate neighbourhood is
Lazarevskoe Cemetery, where Liz saw the graves of Dostoevsky, Tchaikovsky,
Rimsky-Korsakov and other great men, and where Dick was photographed
together with his host Sergey Khrushchev in front of Euler's grave.

Sergey Khrushchev%
\footnoteA{\url{https://official.satbayev.university/en/teachers/khrushchev-sergey},\\
\url{https://mathscinet.ams.org/mathscinet/MRAuthorID/198407}}
works in complex, functional and harmonic analysis. His work
on orthogonal polynomials started after 2000, while based outside Russia.

The Steklov Institute\footnoteA{\url{https://en.wikipedia.org/wiki/St._Petersburg_Department_of_Steklov_Mathematical_Institute_of_the_Russian_Academy_of_Sciences}},
where Dick went for lectures and discussions,
was halfway at
Nevsky Prospekt. The Hermitage, which they visited several times, was at the
other end of this avenue. They also visited the Russian Museum and they
made a trip outside the city to the Pavlovsk Palace and to the Catherine
Palace in Pushkin.

Liz was struck by the difference between their treatment, as privileged
guests of the Academy, and the daily life of the common people, who were
standing in line all the time and were served mean food.
She observed the surveillance in the hotel by women sitting at desks near the
elevators on each floor and taking your room key when you leave the floor and
returning it when you come back.

\paragraph{September 10--21: Stay in Moscow}\qquad\sLP
Their Hotel Akademicheskaya is near Leninskiy Prospekt, a major avenue running
south-west from the Lenin Monument at the Garden Ring to the (outer) Ring Road.
The Red Square and Kremlin are a 45 minutes walk away; Moscow State University
takes 90 minutes.

Just back from Leningrad they had an appointment at the US embassy (not far
from the \emph{White House}, the office of the Russian government). From there
they walked along Kalinin Prospekt (now called New Arbat Avenue) to some major
Moscow sights:
``the Lenin Library, the Tomb of the Unknown Soldier, around the corner into
Red Square, past Lenin's mausoleum and St.~Basil's Church (a truly magnificent
and striking structure) and through the 19th century building that once housed
hundreds of small shops and is now GUM department store.'' Then a shower
started, so they took the subway back to the hotel.

On a later day Liz brought a pleasant visit to the Central Children's Library,
where she handed gifts from Madison: books and cassette tapes. She was
shown around and had plenty of conversation with the Library staff.

In their hotel room Dick and Liz were harassed by bed bugs.
In one night Dick had 13 bites. Liz had less, but she reacted vehemently to
them and she had to consult a doctor at the American embassy.
\mLP
\textsl{And bed bugs it is! Dick was up at 2:30 a.m.\ with six new bites and he
wasn't about to crawl back in there with the beasties who were making a meal
out of him. Can't blame him. However he just left to give his presentation at
the Gelfand seminar and that's a strenuous evening ahead for him with only
three hours of sleep last night. Since talking mathematics is what he's here
for, that may be enough to keep him going.}
\mLP
\textbf{\large September 14, a busy day for Dick}
\sLP
On this Monday Dick, while having slept so badly, was scheduled to
meet successively the three research groups which were relevant for him.
\mLP
\textbf{1. Nikiforov, Uvarov and Suslov}\sLP
According to Liz' diary he was picked up in the morning by
``a trio of physicists and a translator'' and ``one of the physicists
brought me a bouquet of flowers from his garden''.
Actually two of these
(mathematical) physicists were
Arnold Nikiforov\footnoteA{\url{https://tinyurl.com/27vdwjws}
(Russian Wikipedia page),\\
\url{https://www.mathnet.ru/php/person.phtml?personid=27523&option_lang=eng}}
(1930--2005)
and Vasily Uvarov\footnoteA{\url{https://tinyurl.com/yuramt9m}
(Russian Wikipedia page),\\
\url{https://www.mathnet.ru/php/person.phtml?personid=24913&option_lang=eng}}
(1929--1997), both working at the Keldysh Institute.
Nikiforov brought the
flowers and his daughter acted as translator. In front of the Kurchatov
Institute\footnoteA{\url{https://en.wikipedia.org/wiki/Kurchatov_Institute}}
Sergei Suslov\footnoteA{\url{https://search.asu.edu/profile/90437},\\
\url{https://mathscinet.ams.org/mathscinet/MRAuthorID==195261}} (1955)
also joined them.
As a security measure, only employees were allowed to enter that building.
Therefore the meeting took place in a nearby international office building.
As described in the previous section, 
these three men were very active in discrete classical orthogonal
polynomials, culminating in the Russian edition of
the book \cite{N-S-U85b}, so there was a lot to discuss.

Dick met Natig Atakishiyev,
a frequent coauthor of Sergei Suslov, at the end of his Moscow trip.
Only on that occasion Dick told about Jim Wilson's biorthogonal
rational functions \cite[Chapter~III]{Wilson78}/
\mPP
\includegraphics[height=4cm]{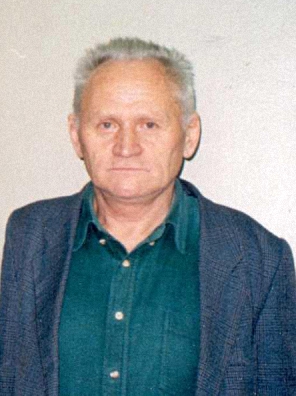}
\hskip2cm
\includegraphics[height=4cm]{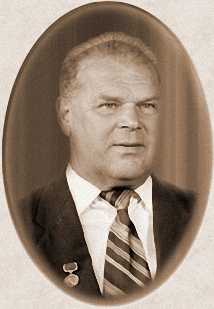}
\hskip2cm
\includegraphics[height=3cm]{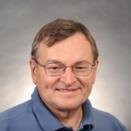}

{\small \hskip1cm Nikiforov\hskip3.5cm Uvarov\hskip4cm Suslov}
\mLP
\textbf{2. The Steklov Institute (Nikolsky and Gonchar)}\sLP
The next visit was to the prestigious Steklov
Institute where Dick met the group headed by Sergey
\mbox{Nikolsky}\footnoteA{\url{https://en.wikipedia.org/wiki/Sergey_Nikolsky},
\\
\url{https://www.mathnet.ru/php/person.phtml?personid=9182&option_lang=eng}}
\cite{Nikolsky03} (1905--2012)
and Andrey Gonchar\footnoteA{\url{https://en.wikipedia.org/wiki/Andrey_Gonchar_(mathematician)},\\
\url{https://www.mathnet.ru/php/person.phtml?&personid=8518&option_lang=eng}}
(1931--2012).
Sergei Suslov, who had been an undergraduate student of Nikolsky,
was also allowed to join there.
An informal seminar was held with
Dick as the main speaker. Nikolsky presented Dick his two-volume analysis book
in English translation \cite{Nikolsky85},
which had recently been published by Mir Publishers. Many years later,
when Sergei Suslov had settled in the USA, Dick offered Sergei these books,
which Sergei had studied from the Russian original as a student.

Next Gonchar (who would be during 1991--1998 vice president of the Russian
Academy of Sciences) invited the small group of the seminar to his home
for lunch.
\mPP
\hskip2cm
\includegraphics[height=2.5cm]{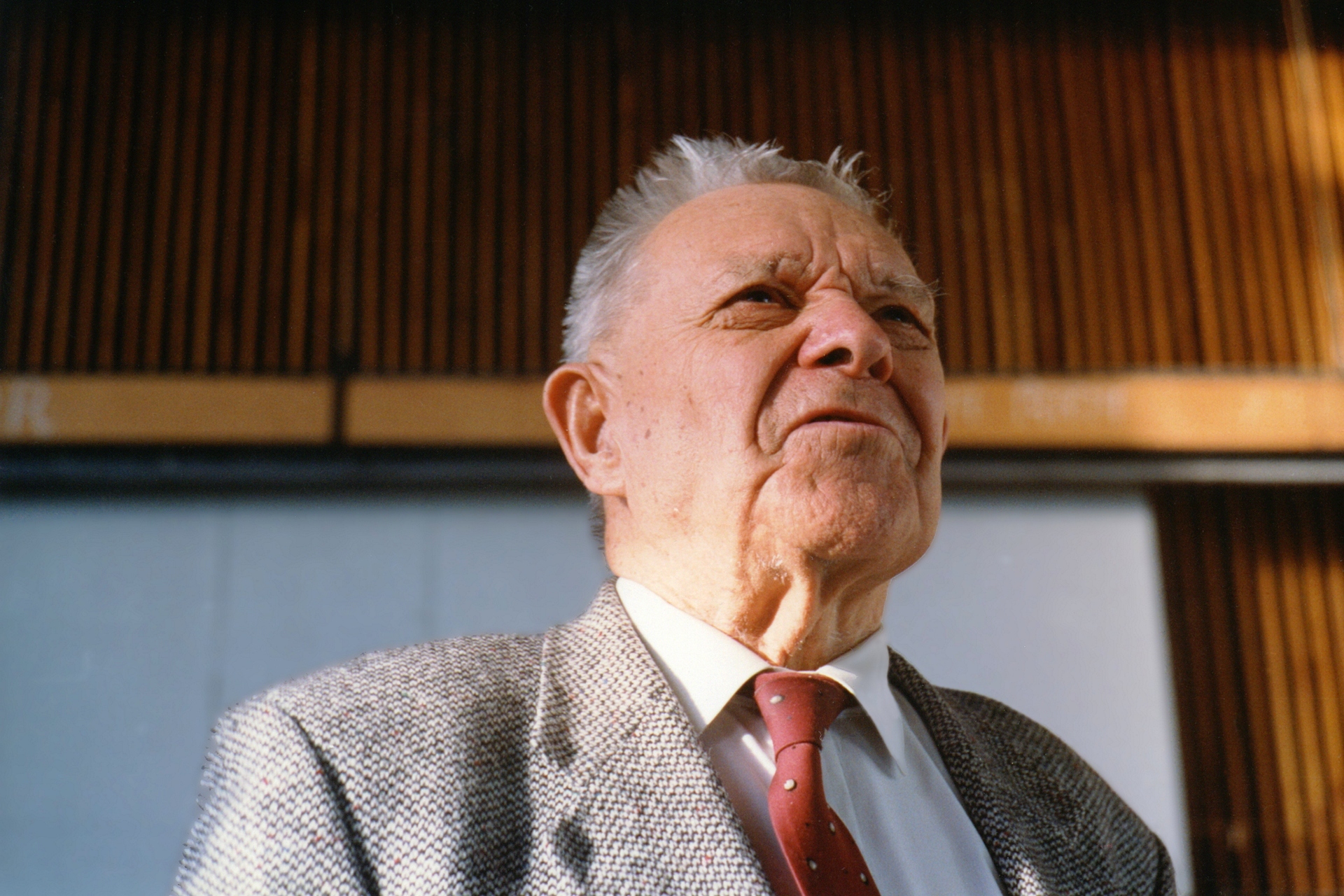}
\hskip2cm
\includegraphics[height=2.5cm]{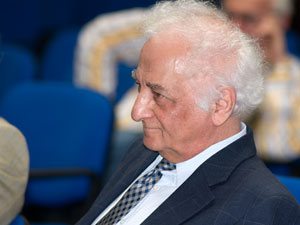}

{\hskip3.7cm\small Nikolsky\hskip4.3cm Gonchar}
\mLP
\textbf{3. The Gelfand seminar}\footnoteA{Read more about the Gelfand
seminar in the Foreword by S. Gindikin to \cite{IMGelfandSeminar} and in
\cite[pp.~25--26, 48, 164]{Retakh13}}
\mLP
\begin{minipage}{9cm}
In the evening Dick spoke in the Gelfand seminar at Moscow State
University (where unfortunately Suslov could not get through the security).
At this university Gel'fand had a group with, among others, Zelevinsky and
Levitan, which had just started to develop Gel'fand's version of generalized
hypergeometric series.
His group did not communicate much with the
group at the Steklov Institute (note also that Gel'fand,
although elected to the
Soviet Academy of Sciences in 1953 as a corresponding member,
did not receive full membership in the academy until
1984 \cite[p.25]{Retakh13}). 
\end{minipage}
\begin{minipage}{7cm}
\includegraphics[width=7cm]{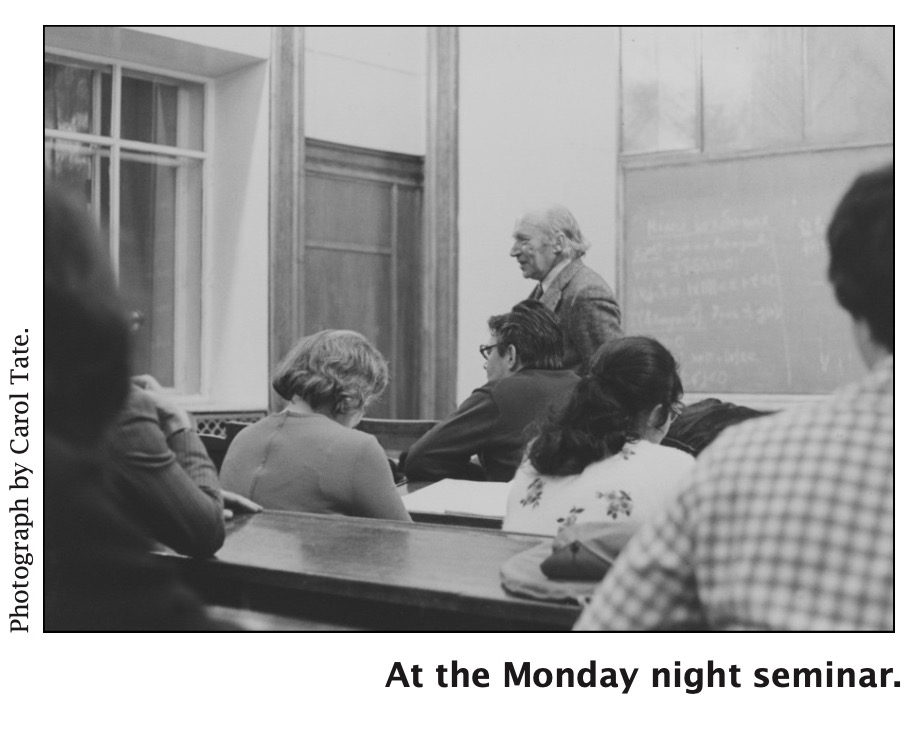}
\end{minipage}
\paragraph{Further activities}
During the stay in Moscow Liz and Dick together or Liz alone visited such
places as Ostankino Palace, Kolomenskoye Museum-Reserve, St.~Basil's 
Church, the Church of St.~George, the Kremlin, Novodavichi Convent and
Nikitniki Cathedral. Liz also brought children's books, taken with her from
Madison, to the Central Children's Library and she was shown around there.

Dick and Liz were not satisfied with the food and the service in the
restaurants, but they had some very good meals in the homes of Russian
colleagues:
\bLP
\begin{minipage}{9cm}
\textsl{The Zelevinskys invited us back to their apartment for a meal,
which was the culinary high point of our stay in Russia. We had borsht, a
mushroom and onion pie, a cheese pie, a potato and meat dish which reminded me
somewhat of scalloped potatoes, a vegetable "sauce" which included eggplant and
tomatoes and was eaten on a slice of bread, fresh sweet watermelon, and, for
dessert, three kinds of cakes -- two made with red currents and one with
strawberries -- and tea. It was absolutely delicious and we enjoyed every
bite.}
\end{minipage}
\qquad\qquad
\begin{minipage}{4cm}
\includegraphics[width=4cm]{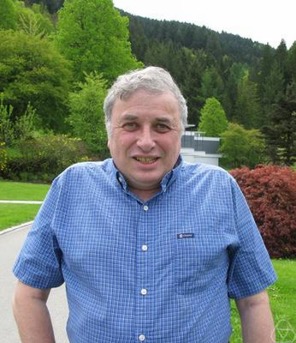}
\end{minipage}
\bLP
Zelevinsky%
\footnoteA{\url{https://en.wikipedia.org/wiki/Andrei_Zelevinsky},\\
\url{https://mathscinet.ams.org/mathscinet/MRAuthorID/191850}}
 (1953--2013) did important work in algebra, combinatorics and
representation theory. He was very much involved in the project on
Gelfand hypergeometric functions which had just started. He emigrated to
the United States in 1990.
\bLP
\begin{minipage}{4cm}
\includegraphics[width=4cm]{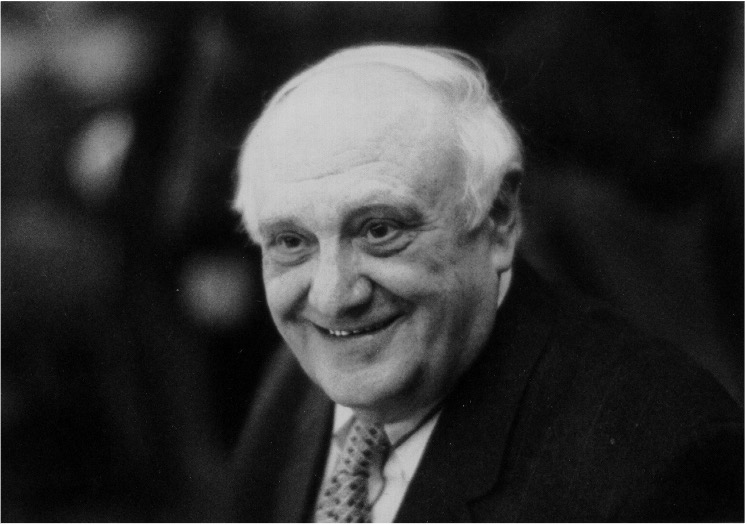}
\end{minipage}
\qquad\qquad
\begin{minipage}{9cm}
\textsl{We had dinner with the Levitans, their daughter-in-law and Zelevinsky
last night. It was a lavish meal with desserts that never stopped, served in
front of a window high over the Moskva River with a marvelous view of Moscow
spread out on the other side of the river.}
\end{minipage}
\bLP
Levitan%
\footnoteA{\url{https://en.wikipedia.org/wiki/Boris_Levitan},\\
\url{https://mathscinet.ams.org/mathscinet/MRAuthorID/198925}}
 (1914--2004) worked on almost periodic functions, Sturm--Liouville
operators and inverse scattering. His paper \cite{GelfandLevitan55} (1955)
with Gel'fand had big impact.
In 1992 he emigrated to the United States.
\mPP
In her diary Liz is aware of Gorbachev's attempts to change the economic
system, but she does not mention \emph{perestroika} and there is no forefeeling
of the big events in 1989--1991 which would tear down the Iron Curtain
and end the
Soviet Union. In general she does not like the Soviet system, but she likes
the Russian culture and the people.
\section{Japan, October 6 -- November 14, 1987}
Schedule:
\begin{itemize}
\item
October 6--12: Sendai (host: Satoru Igari, Tohoku University)
\item
October 12--18: Tokyo (hosts: Suki Kubo, private, and Reiji Takahashi, Sophia
University)
\item
October 18--23: Nagoya (hosts: Kazuhiko Aomoto, Nagoya University, and Toshihiro
Watanabe, Gifu University)
\item
October 23--26: Nagano (conference and excursions with the Aomotos)
\item
October 26--28: Nagoya
\item
October 28 -- November 12: Kyoto (host: Michio Jimbo, Kyoto University)
\item
November 12--14: Tokyo (host: Suki Kubo)
\end{itemize}
\section{Australia, November 15 -- December 10, 1987}
Schedule:
\begin{itemize}
\item
November 15--26: Sydney (hosts: Alfred van der Poorten, Macquarie University,
and David Hunt, University of New South Wales)
\item
November 26 -- December 10: Canberra (host: Rodney Baxter, Australian National
University)
\end{itemize}
\section{India, December 11, 1987 -- January 8, 1988}
Schedule:
\begin{itemize}
\item
December 11--14: Madras (visit to Mrs.~Ramanujan)
\item
December 15--16: Chidambaram (conference at Annamalai University)
\item
December 17: excursion to Kumbakonam (visit to Ramanujan's home)
\item
December 18: Pondicherry (closing session of conference at Annamalai University)
\item
December 19--26: Madras (conference at Anna University and International Conference
opened by Prime Minister Rajiv Gandhi)
\item
December 27--28: Poona (conference)
\item
December 29 -- January 8: Bombay (visit to Tata Research Institute and conference)
\item
January 9--15: London, UK (culture and shopping, no math)
\end{itemize}
Conferences attended:
\begin{itemize}
\item
December 15--18: Ramanujan Centennial International Conference
\cite{ProcAnnamalai88}, Annamalai Nagar, Chidambaram.
\item
December 19--21: International Ramanujan Centenary Conference held at Anna
University, Madras. The conference included a Number Theory symposium
\cite{NumberTheoryMadras89}.
\item
December 22--27: International Ramanujan Birth Centenary Conference,
Institute of Mathematical Sciences, Madras (with interruption on December 24--25
by political turmoil).
\item
December 26--28: Ramanujan Birth Centenary Year International Symposium on
Analysis \cite{RamanujanSymp89}, University of Pune, Poona.
\item
January 4--11: Ramanujan Birth Centenary International Colloquium
\cite{NumberTheoryBombay89}, Tata Institute of Fundamental Research, Bombay.
\end{itemize}
\paragraph{Acknowledgment}
Writing this paper would have been impossible without all information
from personal memory which Sergei Suslov has provided to me. He was an
eyewitness and an active player in Moscow during the 1980s, and his meeting
with Askey in Moscow on September 14, 1987 was crucial for his further
mathematical development and offered to Askey one of the most valuable contacts
during his Russian trip.
Sergei, thank you so much, also for a final reading and commenting
of my draft.

\end{document}